# What the submission does not show: Complexifying cognitive offloading in students' mathematical work with LLMs

Igor' Kontorovich[1] · Miriam N. Wallach[2,3]

[1] Faculty of Education, The University of Hong Kong, Hong Kong SAR, China

[2] Orot Israel College, Elkana, Israel

[3] Technion Israel Institute of Technology, Faculty of Mathematics, Haifa, Israel

**Corresponding author:** Igor' Kontorovich, igork@hku.hk

**ORCID:** Igor' Kontorovich 0000-0002-3353-5445; Miriam N. Wallach 0000-0002-7656-8853

# What the submission does not show: Complexifying cognitive offloading in students' mathematical work with LLMs

## Abstract

A common concern is that using large language models (LLMs) lets students offload mathematical work that tasks are designed to elicit. Yet such claims often rest on outcomes (e.g., submitted solutions), rather than on evidence of students' work with LLMs. We followed 23 first-year linear algebra students as they solved three differently designed homework tasks with LLMs. Drawing on self-regulated learning and choice-affluent environments, we interpreted each of 429 prompts as a trace of a regulatory event and compared 59 chatlogs with the problem solutions the students submitted. Prompts formed five structures, casting the LLM mainly as a solver, explainer, or validator. The mix of structures shifted across the tasks, and ten of the 15 students who completed all three tasks changed their dominant structure. Most final submissions stayed close to LLM output, including those that followed exchanges in which students had explained, conjectured, and checked their own ideas. We conclude that an LLM-aided submission is not a trustworthy marker of students' offloading of the mathematical work, and that offloading needs to be examined in each regulatory space and part by part rather than as a whole.

## 1 Rationale

Consider a homework assignment in a linear algebra course, in which students were encouraged to solve problems using large language models (LLMs). One problem asked students to prove or disprove whether any square matrix can be decomposed into a sum of an upper triangular and a lower triangular matrix. Now consider Figures 1a and 1b: Hanna's solution from an exchange with an LLM and the solution she eventually submitted. The two solutions overlap substantially, with only minor differences. Such submissions were common in the course we studied, as most submitted solutions stayed close to what students obtained from LLMs.

Hanna's submission echoes emerging findings from the mathematics education literature. Noster et al. (2024) found that prospective teachers who solved tasks with ChatGPT often submitted its answers, including incorrect ones. Such submissions fuel a concern that runs through educators' responses to *generative artificial intelligence* (genAI): students offload to LLMs the mathematical work teachers assign. The broader literature captures such practices through the construct of "cognitive offloading," the use of external actions or tools to reduce the cognitive demands of a task (Risko & Gilbert, 2016).

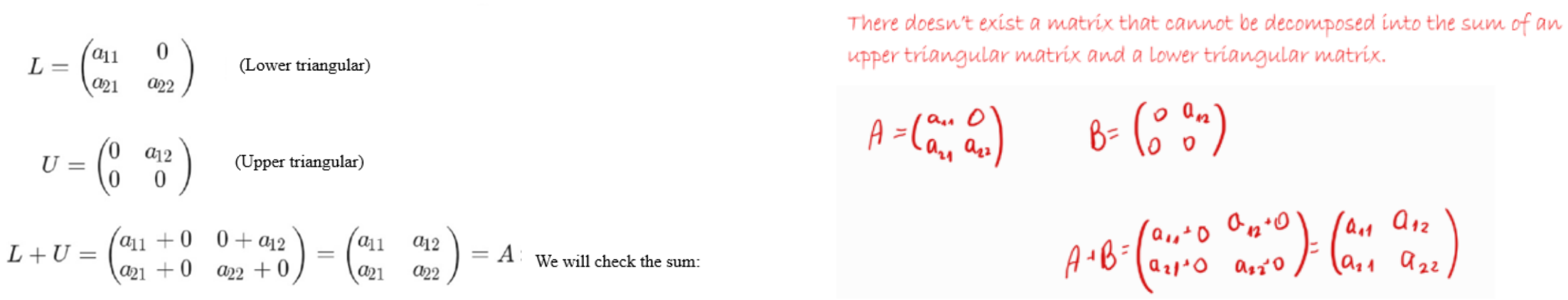

**Fig. 1** Hanna's solution from an LLM (a) and her submitted solution (b). Translated from Hebrew

Our study seeks to complexify this picture. Figure 2 presents an abridged excerpt from Hanna's exchange with the LLM. Before she adopted its solution, Hanna asked the LLM to explain a term from the problem formulation (see [1]), then monitored the obtained explanation to identify another unclear term (see [2]). Hanna also formulated a conjecture and asked the LLM to verify it (see [4]). None of this mathematical work is visible in her submission. Accordingly, our paper seeks to show that a solution developed with an LLM's aid is not a trustworthy marker of the mathematical work that the student invested beforehand. We pursue this aim by examining students' activity within the two spaces: the exchange with the LLM and the transition from that exchange to the eventually submitted solution (see Section 3.3). Indeed, inferences about cognitive offloading commonly rest on outcomes, which speak, at most, to the second space (see Section 2.2).

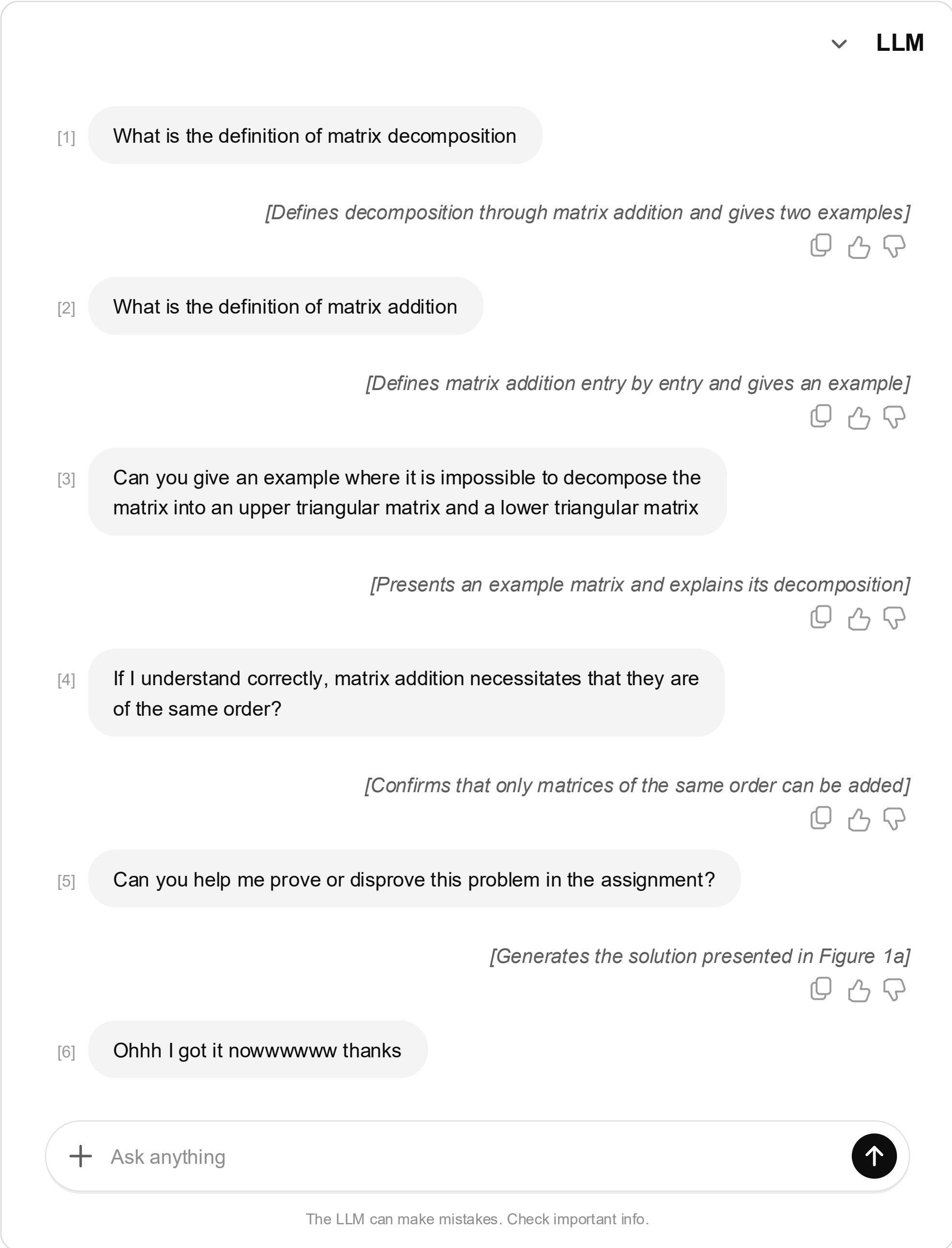

**Fig. 2** An abridged excerpt from Hanna's exchange with an LLM. Translated from Hebrew. Numbers in square brackets mark Hanna's prompts; italic text in square brackets sums up the LLM's replies

Furthermore, we seek to contribute to an emerging line of research that captures the diversity of ways students work with LLMs on mathematical tasks. Studies have examined students and teachers proving with LLMs (Yoon et al., 2024; Zhuang, 2026), modeling with them (Borba et al., 2026; Cevikbas & Kaiser, 2026), and evaluating their mathematical products (Aydoğdu & Girit Yıldız, 2026; Caylan-Ergene et al., 2026; Noster et al., 2024). However, most of this work was arranged by researchers in single sessions, which may not mirror how students use LLMs across their courses (see Section 2.1).

We approach students' mathematical work with LLMs through self-regulated learning and choice-affluent environments (Koichu, 2018), which let us read prompting as regulation within the learning paths that a course extends once an LLM becomes a legitimate knowledge resource (Section 3). Our study took place in a linear algebra course where students regularly used LLMs for homework, which let us follow the same students across differently designed tasks and compare their exchanges with LLMs against the solutions they submitted.

# 2 Research background

## 2.1 Student–LLM interactions and the settings in which they have been studied

Mathematics education came to genAI through the question of what the technology can do (Pepin et al., 2025a), and only afterward through the question of what students do with it (Pepin et al., 2025b). Empirical accounts of the interaction itself nevertheless appeared almost at once. Noster et al. (2024) collected chatlogs from eleven pre-service teachers who solved four mathematical tasks with ChatGPT and coded the "prompting techniques" they employed. The study noted a tendency toward "zero-shot" prompting, in which participants entered the assigned task and accepted the response. Yoon et al. (2024) interviewed three undergraduates proving six mathematical statements and built a framework to explain students' work with LLMs. Dilling and Herrmann (2024) read pre-service primary teachers' geometry work through instrumental genesis and identified prompt patterns at three levels, from the single prompt to the chat interaction taken whole. Kock et al. (2025) analyzed engineering students' accounts of ChatGPT through Vergnaud's notion of "use schemes." Borba et al. (2026) treated ninety-four students in a differential equations course as "collectives of humans-with-GenAI."

Two features of this body of work are worth noting. The first is the length of the engagement. With the exception of Borba et al. (2026), most of them follow a single session in which participants work on tasks set by the researchers. The second is participants' stance. Many were prospective teachers who were asked to inspect the technology as much as to learn with it. In other words, when a pre-service teacher prompts an LLM to solve a mathematical problem, they also evaluate its usefulness for their future teaching (Caylan-Ergene et al., 2026). Yet students increasingly work with LLMs on their own: in a survey of almost 700 Irish students conducted in March 2025, 91% reported using generative AI for mathematical purposes (Ní Fhloinn & Feeney, 2026). A handful of studies have begun to follow LLM use within a course: among upper-secondary students in an algebra course (Brasken et al., 2025), among undergraduates in calculus over a semester (Zhuang, 2025, 2026), and among prospective teachers in a methods module that built AI into its assessment (Mahlaba, 2026). Yet even these studies rarely examined students' work with LLMs across several assignments that differ in design. Our participants, though prospective teachers, took part as mathematics students in a course where LLMs were a regular component of homework, which let us study how their work with LLMs varies from one assignment to the next.

## 2.2 Cognitive offloading and what outcomes can show

Educators' concern regarding genAI is that students can hand their mathematical thinking over to it. Risko and Gilbert (2016) treat offloading as neutral, and often as a sound strategy, since relieving working memory of one demand can free it for another. The concern about genAI is narrower: that it lets students offload precisely the mathematical work that a task was designed to elicit.

Several studies report findings in this direction. Within mathematics education, Yoon et al. (2024) observed students accepting AI-generated proofs without engaging with them critically. Bulut and Borromeo Ferri (2026) found that tenth graders relied almost exclusively on entering the task as given and rarely validated what ChatGPT returned. Similarly, half of the prospective teachers Slobodsky and Durcheva (2025) surveyed did not check LLM-generated solutions. Outside mathematics education, experimental research has supported this concern. In a field experiment with nearly 1,000 high school students, those who practiced with an unrestricted GPT-4 tutor performed better during practice but worse than a control group on a subsequent examination taken without it (Bastani et al., 2025). University students who researched a scientific issue with ChatGPT rather than a search engine reported a lower cognitive load but produced weaker justifications (Stadler et al., 2024). Closer to our setting, Wang et al. (2025) infer “epistemic offloading” in undergraduate linear algebra and calculus from the gap between students’ homework and examination performance. With few exceptions, however, researchers infer offloading from outcomes such as submitted solutions, test scores, or survey responses, which speak, at most, to students’ transition from working with LLMs to the submitted solutions, and not to the preceding work.

### 2.3 The prompt as a unit of analysis

A growing body of work treats students’ prompts as data in their own right. Within mathematics education, several of the studies reviewed in Section 2.1 coded prompts (Dilling & Herrmann, 2024; Noster et al., 2024). Aydoğdu and Girit Yıldız (2026) related the types of prompts that pre-service teachers wrote to the quality of the outdoor mathematics activities they designed. Gülburnu and Oğuz (2026) reported that what tenth graders learned about functions with ChatGPT was tied less to the model’s responses than to the quality of their own questioning and verification.

Empirical analysis of prompts has developed faster outside mathematics education. In computing education, Sheese et al. (2024) categorized more than 2,500 queries submitted to an LLM-powered assistant over a twelve-week deployment. Güner and Er (2025) derived five interaction profiles from students' prompting in a programming course. The first feature of these studies is that their categories tend to describe how a prompt is written, the technique it employs, or the topic it addresses. They are read less often for the mathematical work that the student assigned to the LLM and the work that the student kept. Accordingly, the division of labor that a prompt sets up remains largely undescribed, and it is this division that our study aims to capture (see Section 4.3).

A second feature concerns where prompting is located. Prompting is often described as a property of students: whether as a tendency of a group of participants (e.g., zero-shot prompting in Noster et al., 2024), or as a profile of individuals (e.g., Güner & Er, 2025). In the single-session settings reviewed in Section 2.1, a way of prompting cannot be told apart from the task that elicited it. It thus remains unclear whether prompting constitutes a personal style or a response to the task one is working on.

A third feature concerns where the analysis stops. Prompt studies examine the exchange with the LLM, and they rarely follow what students do with its products afterward. Noster et al. (2024) are a partial exception, since they noted that participants often submitted the answers that ChatGPT gave them. Outside mathematics education, Sawalha et al. (2024) related students' prompting strategies to the accuracy of the model's responses and to students' performance on a final examination. Yet an examination score, like the outcomes reviewed in Section 2.2, does not show how a product of the exchange entered the work that the student presented as their own. Our study brings the exchange and the outcome together by comparing each chatlog with the solution that the student eventually submitted.

The fourth feature is that the frameworks used in prompt analysis were adapted from other areas. Sheese et al. (2024) drew on theories of academic help seeking (Karabenick & Berger, 2013; Nelson-Le Gall, 1981). Fan et al. (2024) coded interaction traces against Zimmerman's (2000) cyclical model of self-regulated learning. The help-seeking literature, for instance, distinguishes between "executive" help seeking, in which the helper is asked to do the work, and "instrumental" help seeking, in which help is sought so that the learner can do the work themselves (Nelson-Le Gall, 1981). Yet unlike the helpers in the original frameworks, an LLM returns text that can be submitted as it is, which makes the passage from the exchange to the submitted product a site of regulation in its own right. Accordingly, we draw on self-regulated learning, as Fan et al. (2024) did, but we read the content of students' prompts and extend the analysis to the transition from the exchange to the submitted solution (see Section 3.3).

## 3 Theoretical framework

Our theoretical framework coordinates two bodies of work that have developed largely independently: self-regulated learning (Section 3.1) and choice-affluent environments (Section 3.2). Section 3.3 adapts them to our setting.

### 3.1 Self-regulated learning

*Self-regulated learning* (SRL) is an umbrella term for the processes by which learners set goals for their own learning and then govern their thinking and actions to serve those goals. Within mathematics education, De Corte et al. (2011) argued that "self-regulation constitutes a major characteristic of productive mathematics learning" (p. 155).

SRL has origins in social cognitive theory, particularly in Bandura's (1991) triad of self-observation, self-judgment, and self-reaction. The triad operates within a reciprocal relation between person, their actions, and the environment. From this source, SRL inherited the commitment that regulation is contextual, meaning learners can regulate differently across tasks and settings.

Winne and Perry (2000) distinguished between studying self-regulation as an “aptitude” and an “event.” The former focuses on relatively enduring attributes expected to hold across occasions. The event lens centers on regulation as it emerges within a stretch of activity bounded by what preceded it and what follows. Given our focus on students’ mathematical work with LLMs, we adopt the event view throughout.

Among social-cognitive accounts of self-regulation, Zimmerman’s (2000) cyclical model of forethought, performance, and self-reflection is widely used. Winne and Hadwin (1998) offered an information-processing account of a different character. They described four weakly sequenced phases, running from task definition through goal setting and planning to the enactment of tactics and to larger-scale adaptation. The researchers specified each phase as a configuration of conditions, operations, products, evaluations, and standards, so that comparing a product against a standard drives any control the learner exercises.

Boekaerts (1999) widened the models’ scope, arguing that learners pursue a mastery pathway and a well-being pathway. Learners often pursue both at once, and regulation then consists of managing the tension between the two. In other words, analysts should not presume that students’ regulation is learning-oriented by default. For instance, a learner can disengage from an activity to protect their sense of competence, which still instantiates regulation. Accordingly, the goal against which a learner judges a product is not always the goal the task was set to serve.

Panadero’s (2017) comparison of SRL models identified four recurrent constructs. *Task definition* supplies the construal the learner forms of what is to be done. *Standards* provide the criteria against which a product is to be judged, and they are ordinarily tacit. *Monitoring* is understood as comparing a product against a standard. The *recursion of products into conditions* accounts for how what one phase generates enters the next as part of the situation to which the learner must respond.

### 3.2 Choice-affluent environments

Koichu (2018) proposed a model of "shifts and choices" as a confluence of research understanding of how mathematical problem solving occurs and practical considerations for promoting it. The model centers on finding solutions to challenging problems within a configuration of resources available to the solver. A foundational premise of the model posits that a solver's pathway to a problem solution

> is stipulated by choices the solver is empowered to make and by enacting the following types of resources: (i) individual resources, (ii) interaction with peer solvers [...], and (iii) interaction with a source of knowledge about the solution or its parts, such as a textbook, an internet resource, a teacher, or a classmate [...] (pp. 310–311).

Nearly a decade ago, Koichu (2018) argued that these three resource types capture "all frequent situations of problem solving in instructional settings" (p. 311). We extend the third point by positioning LLMs as a contemporary source of knowledge.

Our extension comes with two caveats. First, LLM-generated mathematics is not always valid or didactically apt. Yet the same holds for the other sources in Koichu's (2018) list. Second, Koichu describes a source that acts "as a provider of cues to the solution ... but not as a source of telling the solution" (p. 317). An LLM, in contrast, will tell a full solution when asked. We still see it as fitting the model, since what an LLM discloses depends on what the student prompts for. Hence, how much the LLM discloses becomes a choice the student makes when working with it.

An environment is "choice-affluent" to the extent that it empowers solvers to make informed choices among the conditions under which they will work and to revise those choices as their work proceeds. Affluence, on this account, is a property of what is extended to the solver rather than of what any solver takes up.

Koichu (2018) defined affluence by five choices. The first concerns the challenge to be pursued, and it ranges not only over which problem to attempt but over the grain at which to attempt it: "from solving a difficult problem to comprehending a worked-out example" (p. 320). The second choice concerns the mathematical toolkit for dealing with the challenge. The third choice concerns the extent of solvers' collaboration, which can range from joint work to solving independently.

The last two choices are especially relevant to our study. The fourth choice concerns the solvers' mode of interaction, that is, "whether to talk, listen, or be temporarily disengaged from the collective discourse, as well as whether to be a proposer of an idea, a responder … or a silent observer" (Koichu, 2018, p. 320). The fifth choice concerns the agent to learn from, which Koichu glosses as "the opportunity to choose whose and which ideas are worthwhile of their attention" (p. 320).

Let us conclude this section with three notes, two of which are methodological. First, Koichu (2018) positioned his model as exploratory, supporting it with single empirical cases. He also highlighted that the model's research usefulness requires methodologies that can trace students' thinking (or shifts of attention, in his terms) across different problem-solving contexts. Accordingly, the next section proposes a methodology for problem solving with LLMs. Second, the five choices are not dimensions along which an environment is to be scored. We therefore take the construct as a preliminary vocabulary for describing what working with LLMs extends to students, keeping it distinct from what they then do with these extensions.

The last note concerns the networking between SRL and choice-affluent environments. Using Bikner-Ahsbahs and Prediger's (2014) terminology, we aim to coordinate the two frameworks. This strategy fits frameworks that offer complementary views of one phenomenon and whose core principles are compatible. Our phenomenon is a student's work on a mathematical problem with an LLM. On principles, the two frameworks agree: SRL posits a reciprocal relation between the learner and the environment (Bandura, 1991), and Koichu's (2018) model places a solver's pathway within the resources an environment makes available. Both also attribute agency to a person to make choices. The frameworks complement the foci of each other. The choice-affluence lens describes what an environment extends to solvers and leaves open what they do with it. The event view of SRL describes what learners do, as they set goals and judge products against standards. Together, the two lenses let us ask which choices working with an LLM makes available and how a student regulates them. We do not merge the frameworks into one theory. Rather, each answers a different part of our questions.

### 3.3 Mathematical work with LLMs and research questions

We adopt a regulatory lens to study students' mathematical work with LLMs, treating an episode of such work as a regulatory event in the sense of Section 3.1.

We situate students' work with LLMs within the broader context of their course studies. Indeed, students can turn to LLMs as a response to a specific task set up by a course they take or as part of their preparation for something specific related to their studies. As such, students choose specific tasks, goals, and strategies for their work with LLMs, while monitoring how this work brings them closer to achieving their goals. The broader context may also specify standards of use, for instance through course instructions about what may and may not be put to an LLM. Within this context, we identify two regulatory spaces: students' exchanges with the LLM and the transition from that exchange to the product the student decided to submit.

The first regulatory space concerns student work with LLMs, focusing on the exchanges between students' prompts and LLM-generated products. Students can choose the mathematics at the focus of the exchange (cf. challenge in Koichu, 2018), articulate the task for the LLM, specify the mathematical toolkit to use or avoid, and even set the standards against which the products of the exchange will be judged. Through prompts, students can also decide how to distribute roles between themselves and the LLM, setting the degree and mode of collaboration.

Students can also decide to direct the LLM to use specific mathematical agents as a basis for its products, such as certain mathematical websites or textbooks (cf. Koichu, 2018). Students regulate by judging LLM products using previously mentioned or new standards, which may relate to the broader course context. These approaches can revise goals, standards, and modes of collaboration, creating a recursion in which students translate judgments of LLM products into new conditions.

The second regulatory space arises when the work with an LLM is expected to inform some final mathematical product that students need to submit. There, students need to decide how to adopt, adapt, or abandon the products and aspects of their previous work with LLMs. Once again, in this transition, students need to choose their goals, which may differ at this point from those that initially shaped the activity and their exchanges with LLMs. Furthermore, more often than not, the broader course context will position the student as ultimately responsible for the final product, with its quality potentially affecting their course studies (e.g., a mark). Hence, students may apply different standards when using the products of their exchanges with LLMs to prepare the eventual product with which they will be associated.

Building on the above, this research aims to explore undergraduate students' mathematical work with LLMs on homework assignments through the lens of regulation, treating a course that legitimizes LLM use as a choice-affluent environment. Specifically, we aim to characterize students' regulation and choices as they work on mathematical problems with LLMs.[1] We operationalize this aim into two sets of questions that correspond to the two regulatory spaces:

(i) What mathematics, toolkits, agents, tasks, and standards can be discerned from students' prompts to LLMs, and what roles do these prompts assign to the LLM and reserve for the student? Also, how do these prompts combine into structures within and across chatlogs and tasks?

[1] The "self" has been questioned in SRL, since regulation often unfolds with and through others (Hadwin et al., 2018). In our case, the LLM is the other. Hence, we use "regulation" rather than "self-regulation."

(ii) How do students adopt, adapt, or abandon the products of their exchanges with LLMs when preparing the problem solutions for submission? How does this uptake relate to students' prompting?

# 4 Methods

## 4.1 Setting and participants

Our data came from a first-year linear algebra course at a teacher education college in Israel. The course covered standard introductory content on systems of linear equations, matrix algebra, and vector spaces. The teacher embedded genAI in the course's regular teaching and learning practices and encouraged students to consult it when they struggled with the course material. Also, genAI played a visible role in homework assignments. Each assignment consisted of several problems, one of which was explicitly designated to be solved with LLMs. The college had no license for a specific LLM, so students could choose their preferred platform. Most students used ChatGPT, while three and five students used Claude and Gemini, respectively. Accordingly, we use "LLM" throughout the text.

The course students majored in secondary mathematics education. Participation in the research was voluntary and separate from the course requirements. All students completed the assignments as part of their coursework, but we included only the submissions of 23 who gave informed consent in the data corpus.

## 4.2 The tasks

The three LLM-involving tasks were distributed throughout the semester. Two design principles were common to all of them. First, the instructions framed the LLM as a learning aid rather than a "homework solver," in line with the view of AI as a learning collaborator (Ouyang &

Jiao, 2021). Second, in addition to their final solutions, the students were instructed to submit screenshots of their complete dialogues with the LLM. Participation varied across the assignments: 19, 21, and 19 of the consenting students submitted chatlogs for Tasks 1 to 3, respectively, since some students did not submit an assignment.

The mathematical content and specific guidelines for using LLMs varied across tasks (see Table 1). Task 1 centered on zero divisors, asking students to prove or disprove whether non-zero 3×2 and 2×3 matrices exist whose product is the zero matrix. Students were asked to lead the mathematical work and received minimal guidance for their exchanges with LLMs. Task 2 focused on decomposing matrices into sums of special types, such as upper- and lower-triangular or symmetric and skew-symmetric matrices. This task provided students with a pool of prompts that the literature considers conducive to learning. Task 3 shifted the focus from producing a solution to comprehending a given one. Students received a step-by-step proof that AB is invertible if and only if A and B are invertible. They were asked to use an LLM as a dialogic partner to explain why each step is valid and how it advances the proof.

**Table 1** Mathematical tasks and guidelines for LLM use

| | Task 1 | Task 2 | Task 3 |
|---|---|---|---|
| Topic | Zero divisors | Decomposing matrices into sums of special types | Invertibility of a matrix product |
| Mathematical task | Prove or disprove: there exist non-zero 3×2 and 2×3 matrices whose product is the zero matrix | Prove or disprove: any square matrix can be decomposed into the sum of (a) an upper and a lower triangular matrix; (b) a symmetric and a skew-symmetric matrix; (c) a triangular and a symmetric matrix | Explain each step of a given proof that AB is invertible if and only if A and B are invertible |
| Guidelines for LLM use | Lead the mathematical work; | A pool of learning-promoting prompts (e.g., “Could you break it down into | Use the LLM as a dialogic partner to explain why each |

| | Task 1 | Task 2 | Task 3 |
|---|---|---|---|
| | minimal guidance for the exchange | simpler subproblems?") | step is valid and how it advances the proof |

### 4.3 Data corpus and analysis

Our data corpus consisted of screenshots of student-LLM chatlogs and the final solutions students submitted as part of their assignments. Twenty-three students contributed at least one chatlog, and 15 of them did so for all three tasks. The corpus comprised 59 chatlogs containing 429 mathematical student prompts, which amounts to 7.3 prompts per chatlog on average. We counted as mathematical all prompts that dealt with mathematics and excluded those that only greeted or thanked the LLM. The chatlogs and the solutions were in Hebrew. We translated the excerpts we present for illustrative purposes. All student names are pseudonyms.

In the first phase of the analysis, we examined each mathematical prompt through the lens of the constructs in the first research question. We also noted whether each chatlog centered on the assigned problem or on mathematics that the student set. Every prompt specified the mathematics and the task for the LLM, while other elements were rare (see Table 2). Hence, in the second phase, we grouped recurring combinations of the prompt elements into what we termed *prompting structures*.

We drew the term from Goldin et al.'s (2011) engagement structures, defined as "an idealization involving a characteristic motivating desire or goal, actions including social behaviors toward fulfilling the desire [...] and possible outcomes" (p. 547). Borrowing Goldin et al.'s spirit and architecture of the construct, we define a "prompting structure" as an idealized constellation of an LLM's task that can be discerned from a student's prompt and the mathematics the prompt specifies. Following Goldin et al. (2011), we name the structures in students' voice (see Table 3).

A note on what prompts can and cannot show is in order. We interpret each prompt as a trace of a regulatory decision rather than as the decision itself. A prompt records the task that the student set for the LLM, the mathematics they put to it, and the roles that their request assigned, and these are the elements on which our coding rests. However, a prompt does not record the goal that the student pursued, the standards against which they judged the LLM's previous response, or the reasoning that took place between one prompt and the next. For instance, a request to explain a term may arise from monitoring an unclear response, from compliance with the task instructions, or from an attempt to obtain text for submission, and the prompt alone does not decide among these readings. Accordingly, when we refer to students' goals, judgments, or standards in what follows, we refer to what their prompts are consistent with, and we treat such readings as interpretations that the chatlogs support rather than as observations of regulation itself.

Below, we illustrate the structure analysis with Hanna's chatlog from Task 2 (see Figure 2). Prompts 1 and 2 asked the LLM to define a term, first from the problem ("matrix decomposition") and then from its reply ("matrix addition"), so we coded both as "Explain this term." In Prompt 4, Hanna asked the LLM to judge her conjecture about the order of matrices, which we coded as "Is my idea correct?" In Prompt 5, she handed the problem to the LLM; hence "Solve this problem." Prompt 3 was harder to code. It could seek to clarify what decomposition means, or it could ask for a counterexample. Since a matrix that cannot be decomposed would disprove the claim, we coded it as "Solve this problem." Prompt 6 only thanked the LLM, so we did not count it as mathematical.

At the chatlog level, we recorded which structures each chatlog used at least once, how many it combined, and its length in prompts (see Table 4). We also computed the share of each chatlog's prompts that each structure took up, and we refer to the structure with the largest share as the chatlog's dominant structure. We then assigned each chatlog to one of six types by its mix of structures, from mostly doing to mostly checking. A chatlog was led by a single structure when that structure accounted for most of its prompts and the others appeared only in passing. We considered a chatlog as combining structures when several structures took up a substantial share of its prompts (e.g., Hanna's chatlog). We made this coding qualitatively, considering the sequence of the prompts as well as their shares, and we discussed borderline cases until we agreed. Three types were led by a single structure: "Solve this problem" (17 chatlogs), "Explain this term" (19), and "Is my idea correct?" (6). The other three combined structures: Solve + Explain (7), Explain + Check (5), and Solve + Explain + Check (5). Accordingly, when we refer to chatlogs led by explaining or checking, we mean the 30 chatlogs of the "Explain this term," Explain + Check, and "Is my idea correct?" types, and when we refer to chatlogs that involved checking, we mean the 16 chatlogs of the last three types in this order.

In the last phase, we compared each chatlog with the solution the student submitted, noting where it matched the LLM's responses and where it departed from them. Each pair fell into one of four categories, which refine the choice to adopt, adapt, or abandon the LLM's products framed in Section 3.3. In *copy-paste*, the student submitted the LLM's response unchanged, as a screenshot or a copy. In *minor revisions*, the student rephrased the response, adjusted it to the task, or fixed math that the right-to-left Hebrew layout had scrambled. In these cases, the submission overlapped with the LLM's responses to a large extent. In *major revisions*, the solution differed from the response in its mathematics. For instance, the student used the same procedure to reach a different answer, or a different procedure to reach the same one. Finally, in *abandonment*, the student put the LLM's responses aside, and we found no direct link between them and the submitted solution.

To describe how the structures related to each other and to chatlog length, we used Spearman's rank correlations. Since prompts and chatlogs were nested within students and several cells were small, we did not conduct inferential tests. Accordingly, we report counts, percentages, and correlation coefficients, and we read the patterns as associations within our corpus rather than as estimates of population effects.

The second author led the initial coding and presented the emerging categories, with examples, to the first author. We refined the categories through our discussions. Then, the first author independently coded about 40% of the corpus at each phase. Agreement was 87% ($\kappa$ = .81) for the prompting structures, 92% ($\kappa$ = .90) for the chatlog types, and 85% ($\kappa$ = .70) for the uses of LLM products in the submitted solutions. The lower $\kappa$ in the last phase reflects the predominance of minor revisions (66% of the submissions), which inflates the agreement expected by chance. We resolved disagreements through discussion until we reached consensus.

## 5 Findings

### 5.1 Setting up the math: Shifting from inquiry to the assigned problems

In 55 out of 59 chatlogs, the students centered their exchanges with LLMs on the problem set in the assignment. We discuss these exchanges in detail in Section 5.2. Here, we elaborate on the four exceptions in which the students defined mathematics that was not identical to the one featured in the assigned problems. Figure 3 shows one of them, submitted by a student in Task 1.

Recall that Task 1 asked students to prove or disprove that the product of a 3×2 and a 2×3 matrix can be the zero matrix. The exchange does not address this mathematics. Rather, Maya opened the exchange with a general question about whether any two matrices can be multiplied (Prompt 1). She then turned to products of special matrices, asking whether a diagonal matrix times a general one stays diagonal (Prompt 4) and whether two general matrices can give a diagonal one (Prompt 5). Thus, Maya set her own mathematics and used the LLM to validate or reject her conjectures. True, her prompts related to the focal problem, since the LLM's replies could help her narrow the search for suitable matrices. Yet no prompt mentions zero products or the matrix sizes in Task 1. Moreover, diagonal matrices are square. Thus, we argue that Maya's inquiry yielded insights into products of matrix types that go beyond Task 1.

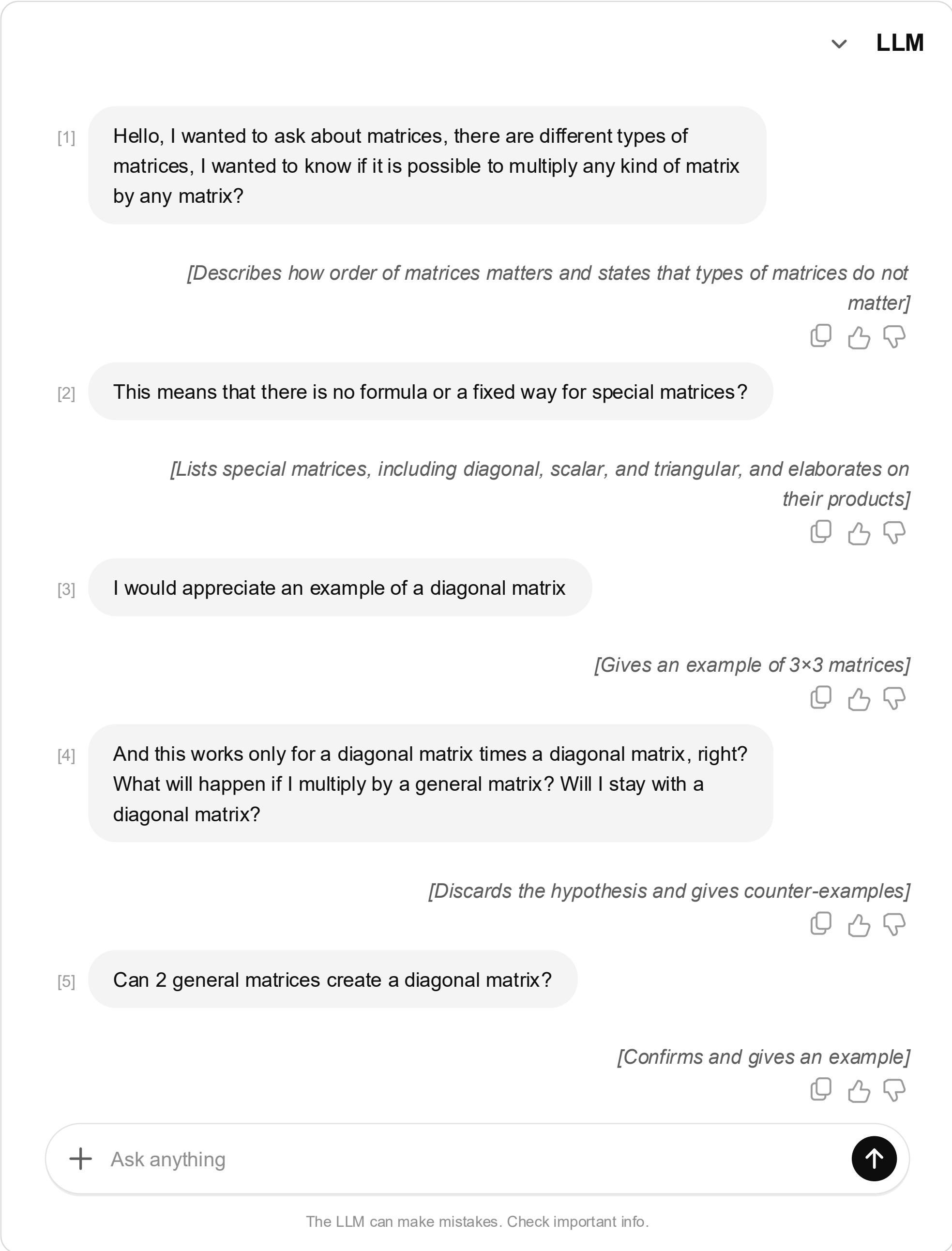

**Fig. 3** Translated excerpt from Maya's chatlog. Numbers in square brackets mark Maya's prompts; italic text in square brackets sums up the LLM's replies

In all four inquiries, the students prompted the LLM with general questions and conjectures that did not mention the specifics of the assigned tasks. Thus, they defined the mathematics themselves, chose which properties to probe, and cast the LLM as the agent to learn from about matrices in general (cf. Section 3.3). The four inquiries also shared a turning point: after a few prompts, the students shifted to the assigned problem, which redefined the task and the roles that the prompts assigned. Accordingly, the assignment did not fix the mathematics the students chose to discuss, since they revised it within the exchange as a regulatory move in its own right.

Notwithstanding, the shift mentioned was not a clean switch from the students' own inquiry to the assigned problem. Even after the turn, the prompts retained some degree of generality that was unnecessary for solving the assigned tasks. For instance, Maya went on to ask whether the product of two general matrices can be the zero matrix and whether a formula exists for constructing such matrices. But she did not mention the 3×2 and 2×3 sizes that Task 1 specified. Her questions thus moved toward the focal problem of zero products while keeping the scope of her earlier inquiry into general matrices. The same held in the other inquiries, where the students approached the assigned mathematics through questions broader than the problem itself. Accordingly, the turn appears as a gradual narrowing, in which the students' self-defined mathematics and the assigned one coexisted for a while. In regulatory terms, the students did not simply give up their own conditions for those set by the tasks. Rather, they adjusted them, possibly in response to LLM-generated products.

### 5.2 Characterizing prompts, prompting structures, and chatlogs

Table 2 summarizes how the constructs from Section 3.3 applied to the students' prompts. They were spread unevenly. The students specified the mathematics at the center of their prompts and stated what they wanted the LLM to do with it. At the same time, no prompt directed the LLM toward particular mathematical resources (e.g., a website, a textbook, course materials). Moreover, students' prompts read as if the LLM produced the mathematics itself. This held even when the LLM-generated responses were unclear or beyond the students' reach, with students referring to the LLM directly (e.g., "I don't understand *your* method," "*your* calculation is wrong").

References to the toolkit were rare (1.4% of the prompts). Only two prompts specified the toolkit in advance, and both named the topic the LLM should draw on rather than a method. The other references emerged in follow-ups, after the LLM used methods unfamiliar to the students. For example, in the exchange in Section 5.1, the LLM built general matrices whose product is the zero matrix using linear dependence of rows. This topic came late in the course, and Maya replied, "I didn't understand this method. What is dependent rows [...]". Thus, students articulated the toolkit in follow-ups that are consistent with judging an LLM product and tracing its misfit to methods beyond their reach. We interpret such instances as traces of regulatory events, in which a student evaluated an LLM-generated product, identified a mismatch, and introduced a new toolkit condition. The condition made explicit a standard that had gone unstated so far.

**Table 2** Elements of the students' prompts

| Element | Mathematics | Toolkit | Agent | Task definition | Standards |
|---|---|---|---|---|---|
| % of prompts | 100 | 1.4 | 0 | 100 | 12 |

*Note.* Percentages are across the three tasks.

Explicit mentions of standards were also infrequent. Only a couple of students consistently asked the LLM to present the math in a "clear" or "simple" way. In some cases, such standards emerged similarly to how the toolkit did, that is, in response to something the student did not understand in the LLM's reply.

We grouped the prompts into five prompting structures, each dividing the mathematical work between the student and the LLM differently (see Table 3). Table 3 shows that the five structures echo the roles set out in Section 3.3. In most prompts, students cast the LLM as a solver, explainer, or validator, and each role left them with a different part of the work. When asked to "solve this problem," the LLM took over the solving, and the student kept only the posing of the task. When asked to explain a term or judge an idea, the LLM supported work that stayed with the student. “Continue my idea” split the solving: the student set the beginning of a solution, and the LLM completed it.

**Table 3** Prompting structures and the roles they assign

| Prompting structure | Description | Roles (LLM / student) | Example |
|---|---|---|---|
| “Solve this problem” | Prompts asking the LLM to carry out a set mathematical task, mostly to solve a problem. | LLM as solver; the student keeps the posing of the task. | “I need you to prove or disprove the following statement: There exists a 2X3 matrix and a 3X2 matrix such that each is not zero but their product is zero. Justify it please” |
| “Explain this term” | Prompts asking the LLM to explain a term, a method, or a step of an argument, whether from the problem, a given proof, or the LLM’s previous reply. | LLM as explainer; the student identifies what is unclear. | “What is the zero matrix?” |
| “Is my idea correct?” | Prompts asking the LLM to judge the validity of the student’s idea, conjecture, or solution. | LLM as validator; the student as solver or author of a conjecture. | “Then if I changed the place of the row of zeroes or the column of zeroes the matrix that comes out wouldn’t be the zero matrix?” |
| “Continue my idea” | Prompts asking the LLM to complete the student’s partial solution. | LLM as co-solver; the student starts the solution. | “For a proof I wrote let A be a square 3 by 3 matrix and let B be a symmetric square 3 by 3 matrix. Give me a hint how to continue from here” |
| “Your math is incorrect” | Prompts rejecting the validity of the mathematics in an LLM response. | Student as validator; the LLM as a solver. | LLM: “[...] decomposition is possible for all square matrices” Student: “But I think that its not true, not every |

| | |
|---|---|
| | matrix can be decomposed” |

Table 4 shows that the mix of structures shifted with the task. “Explain this term” was the most common structure in every task, though only narrowly in Task 1. Requests to check fell, from 33% to 28% and then to 13%. This is consistent with the task designs: in Task 3, the students were given a proof to make sense of, so they had little work of their own to put to the LLM for judgment. Requests to “solve this problem” showed no steady trend. They dipped to 20% in Task 2, where the students had a pool of exemplary prompts, and rose to 37% in Task 3, though from a small base.

Counting chatlogs rather than prompts adds two points. First, “Solve this problem” was more widespread than its share of prompts suggests. It appeared in 84% of Task 1 chatlogs, more than any other structure, and in 70% of all chatlogs. Second, the tasks differed in how many structures a chatlog drew on. In Task 2, nearly every chatlog used “Explain this term” (95%), and the median chatlog combined three structures. In Task 3, most chatlogs were short and used a single structure (58%). Indeed, a student-level view points the same way. Of the 15 students who completed all three tasks, only five kept the same dominant structure throughout: two with “Solve this problem”, two with “Explain this term”, and one with “Is my idea correct?”. The other ten shifted with no common path; for example, from explaining in Tasks 1 and 2 to doing in Task 3. Accordingly, students’ prompting appears to have varied with the task rather than following a fixed personal style (but see Section 6.2).

**Table 4** Distribution of prompting structures across tasks, by prompts and by chatlogs

| Prompting structure | Task 1 | Task 2 | Task 3 | Total |
|---|---|---|---|---|
| *Prompts* | | | | |
| “Solve this problem” | 43 (30.5) | 45 (20.0) | 23 (36.5) | 111 (25.9) |
| “Explain this term” | 49 (34.8) | 105 (46.7) | 32 (50.8) | 186 (43.4) |
| “Is my idea correct?” | 46 (32.6) | 64 (28.4) | 8 (12.7) | 118 (27.5) |
| “Continue my idea” | 2 (1.4) | 8 (3.6) | 0 (0.0) | 10 (2.3) |
| “Your math is incorrect” | 1 (0.7) | 3 (1.3) | 0 (0.0) | 4 (0.9) |
| Total | 141 (100) | 225 (100) | 63 (100) | 429 (100) |

| Prompting structure | Task 1 | Task 2 | Task 3 | Total |
|---|---|---|---|---|
| *Prompts* | | | | |
| *Chatlogs* | | | | |
| “Solve this problem” | 16 (84.2) | 14 (66.7) | 11 (57.9) | 41 (69.5) |
| “Explain this term” | 13 (68.4) | 20 (95.2) | 11 (57.9) | 44 (74.6) |
| “Is my idea correct?” | 13 (68.4) | 16 (76.2) | 6 (31.6) | 35 (59.3) |
| “Continue my idea” | 1 (5.3) | 4 (19.0) | 0 (0.0) | 5 (8.5) |
| “Your math is incorrect” | 1 (5.3) | 3 (14.3) | 0 (0.0) | 4 (6.8) |
| Single structure only | 4 (21.1) | 1 (4.8) | 11 (57.9) | 16 (27.1) |
| Median structures per chatlog | 2 | 3 | 1 | 2 |
| Total chatlogs | 19 | 21 | 19 | 59 |

*Note*. Prompt values are $n$ (%) of prompts within each task. Chatlog values are $n$ (%) of chatlogs in each task that used the structure at least once, so they do not sum to 100. $N$ = 429 prompts in 59 chatlogs from 23 students.

The chatlogs also differed in how the structures combined. The more a chatlog relied on “Solve this problem,” the less it used “Explain this term” (Spearman’s $\rho = -.67$) and “Is my idea correct?” ($\rho = -.58$). In contrast, the use of these two was unrelated ($\rho = -.01$), so they often appeared together. Chatlogs led by “Solve this problem” were also short, with a median of 4.5 prompts, against 7 for chatlogs led by explaining and 8.5 for those led by checking. The 12 chatlogs that used “Solve this problem” alone had a median of two prompts. More broadly, longer chatlogs drew on more structures ($\rho = .61$), and the two rare structures appeared only in long exchanges: all nine chatlogs with “Continue my idea” or “Your math is incorrect” had at least seven prompts. Thus, handing the mathematics to the LLM tended to close the exchange, while explaining and checking kept it going and opened room for other structures, including the rare ones in which students judged the LLM’s output.

## 5.3 From exchanges with LLMs to submitted solutions

We now turn to how students used the products of their exchanges with LLMs in their submitted solutions, along the four categories defined in Section 4.3. Copy-paste submissions, whether rewritten by hand or submitted as snapshots, left no visible traces of the students judging the LLM’s product against their own standards.

Many minor revisions corrected the rendering of the LLM's mathematics, which the right-to-left Hebrew layout had scrambled (see Section 4.3). Such revisions can be read as traces of regulatory events in which students spotted flaws in how the mathematics was rendered, while the solving remained the LLM's. Additionally, the LLM provided long explanations, sometimes using terms not familiar to the students. In such cases, the students selected parts of LLM-generated responses, rephrased and adjusted them to the task. For example, in Task 1, the LLM exemplified zero divisors for Shai with $\begin{pmatrix} 1 & 0 \\ 2 & 0 \\ 3 & 0 \end{pmatrix}$ and $\begin{pmatrix} 0 & 0 & 0 \\ 9 & 8 & 7 \end{pmatrix}$. Shai's submission preserved the matrix $\begin{pmatrix} 1 & 0 \\ 2 & 0 \\ 3 & 0 \end{pmatrix}$ and adjusted the numbers in the second one to become $\begin{pmatrix} 0 & 0 & 0 \\ 4 & 5 & 6 \end{pmatrix}$.

As an example of a major revision, in Task 2, Tamar's LLM argued that given a square matrix $A$, the symmetric $S$ can be built as "$S_{ji} = A_{ij}$ if $i < j$, $S_{ij} = A_{ij}$ if $i > j$, and $S_{ii} = A_{ii}$." The skew-symmetric $T$ can then be defined as "$T_{ij} = A_{ij} - S_{ij}$ if $i < j$, and $T_{ij} = 0$ if $i \geq j$." The solution concluded with "the decomposition is possible for any square matrix, if it is defined in this way." Figure 4 presents the solution that Tamar submitted. There, she borrowed the LLM's approaches to defining matrices in a particular case, which she leveraged to disprove the claim.

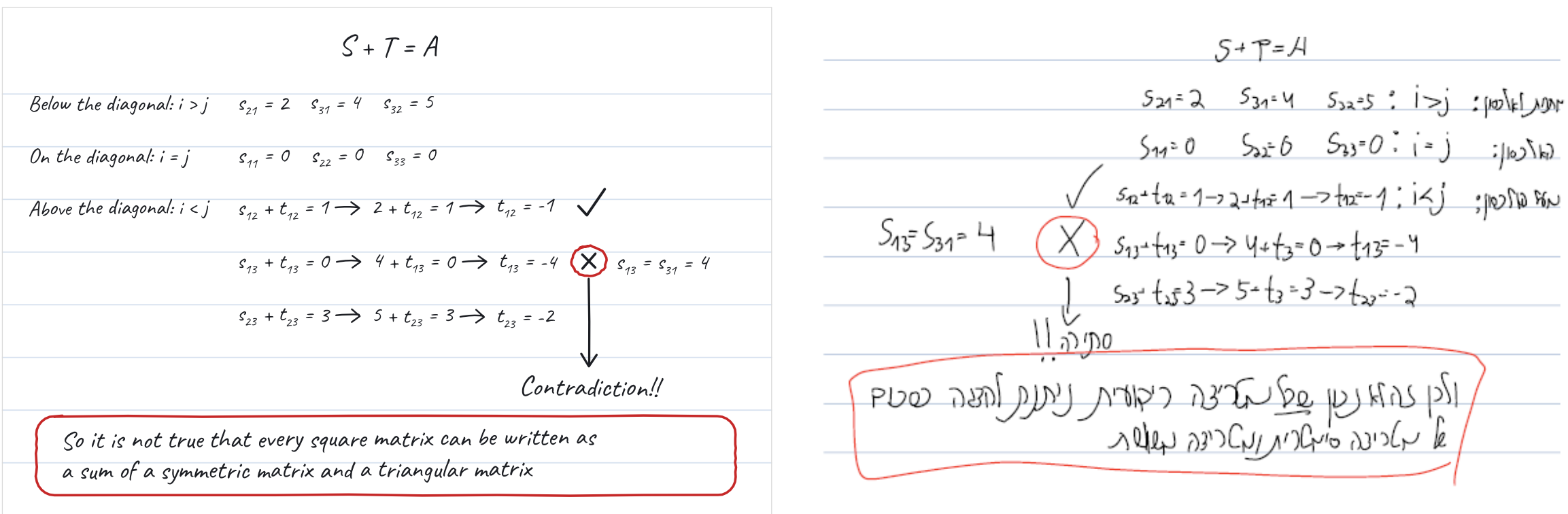

(a) (b)

**Fig. 4** Tamar's submitted solution in Task 2: English translation (a) and original (b)

In one case of abandonment, an LLM-generated solution posited that "it is impossible to decompose a general square matrix into a sum of a symmetric matrix and an upper triangular matrix. The standard decomposition is into symmetric and anti-symmetric matrices." However, the student's solution contained an attempt to prove the claim. Furthermore, several students explicitly stated that their solution discards "AI because it didn't answer correctly." Overall, major revisions and abandonments are clear traces of students monitoring LLM-generated products and deciding that they did not meet their standards.

Table 5 shows that across the collected chatlogs, most submissions drew on what the LLM produced but did not simply reproduce it. In 39 out of 59 cases (66%), the students made minor revisions. In nearly a quarter of cases, the students submitted the response unchanged. Major revisions and abandonment were very rare. Thus, most submissions carried traces of the students judging the form of the LLM's product, yet this judgment seldom reached its mathematics.

**Table 5** Use of LLM products in submitted solutions, by task

| Use | Task 1 | Task 2 | Task 3 | Total |
|---|---|---|---|---|
| Copy-paste | 3 (15.8) | 5 (23.8) | 6 (31.6) | 14 (23.7) |
| Minor revisions | 13 (68.4) | 13 (61.9) | 13 (68.4) | 39 (66.1) |
| Major revisions | 0 | 2 (9.5) | 0 | 2 (3.4) |
| Abandonment | 3 (15.8) | 1 (4.8) | 0 | 4 (6.8) |
| Total | 19 | 21 | 19 | 59 |

*Note*. Values are *n* (%) of submitted solutions within each task.

This profile held across the tasks. Minor revisions accounted for 13 submissions in each task (62–68%). Copy-paste rose from 3 of 19 submissions in Task 1 (16%) to 5 of 21 in Task 2 (24%) and 6 of 19 in Task 3 (32%), although these differences amount to a few cases in each task. Major revisions and abandonment were limited to the first two tasks. In Task 3, where the students made sense of a given proof, no submission departed from the LLM's product beyond minor revisions.

Lastly, we turn to the relationship between students' prompting and their use of LLM-generated responses. Figure 5 shows the main pattern: the more a chatlog drew on explaining and checking, the less often the students submitted the LLM's product as is. Copy-paste made up over half of the "Solve this problem" chatlogs (9 of 17, 53%) but at most one in five chatlogs of any other type, and it was absent from the chatlogs led by checking. Of the 16 chatlogs that involved checking, only one (6%) ended in copy-paste, and 14 ended in minor revisions, including all five that combined the three main structures. Accordingly, students who put their own ideas to the LLM for judgment also tended to rework its products, whereas handing the mathematics to the LLM tended to end in submitting its answer as is. Thus, the roles that students reserved while prompting carried over into how they regulated the move to submission. Taken together, whereas the students' prompting varied across chatlogs and tasks, their submitted solutions varied far less. In 53 of the 59 cases (90%), the submission reproduced the LLM's product or revised it only slightly. In other words, of all the parts of the mathematical work, writing the solution was the one that the students chose to hand over most consistently, keeping for themselves, at most, the editing of its form.

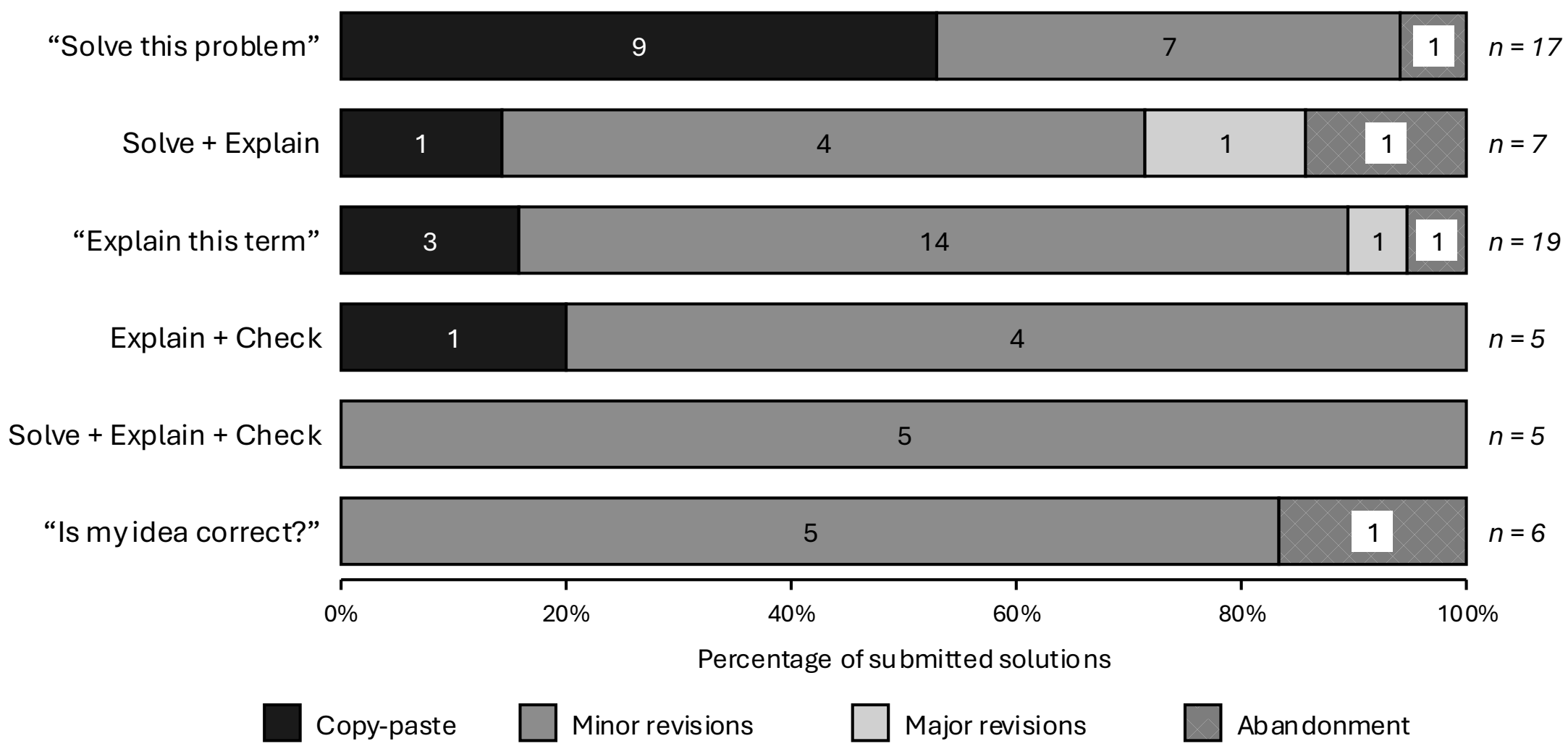

**Fig. 5** Use of LLM products in submitted solutions, by chatlog type. Bars show the percentage of submitted solutions in each chatlog type; numbers inside bars are counts. Chatlog types are listed from mostly doing to mostly checking

## 6 Summary and discussion

This study examined whether a solution developed with an LLM's aid is a trustworthy marker of the mathematical work that preceded it. Our findings show that students divided this work with the LLM in varied ways across prompts and tasks, yet most of their submissions stayed close to the LLM's product.

The main contribution of this study is to complexify the notion of cognitive offloading in the context of students' work with LLMs (Section 6.1). Offloading, as our data suggest, can occur in the exchange with the LLM and in the transition to the submitted solution without the two coinciding. Furthermore, the problematics of this offloading need to be considered in light of the goals that students set for themselves. Two further contributions concern task design (Section 6.2) and what chatlogs record of regulation (Section 6.3). We close with limitations and directions for further research (Section 6.4).

### 6.1 Complexifying cognitive offloading

The concern that students hand their mathematical thinking over to genAI rests largely on outcomes (Section 2.2), which speak, at most, to the second regulatory space. The inference from such outcomes to the exchange itself is precisely what our findings question.

Indeed, the two spaces did not move together. In the 17 chatlogs led by "Solve this problem," the students handed the solving to the LLM, and 16 of their submissions stayed close to its product. The nine of these submitted unchanged most clearly match the concern of cognitive offloading. In the 30 chatlogs led by explaining or checking, most prompts kept the solving with the students, who asked the LLM to explain terms and methods or put their own ideas to it for judgment. Yet 27 of these 30 chatlogs also ended in submissions that stayed close to the LLM's product (see Figure 5). Hanna's case in Section 1 illustrates the point. Before she adopted the LLM's solution, her exchange recorded explaining, clarifying, and conjecturing, and none of this work reached her submission. Thus, a submission that stays close to an LLM's product may follow an exchange in which the student mainly handed the solving over, or one in which the student explained, conjectured, and checked, and none of this work is visible in the submission.

Our findings also suggest that offloading is better described part by part than as a whole. In the exchange, each prompt apportioned the work anew: the student could hand over the solving and keep the posing, hand over the explaining and keep the solving, or hand over the judgment of an idea that was their own. At submission, a similar division appeared. In minor revisions, the students kept the editing, rephrasing the LLM's solution, trimming it, and correcting how its mathematics was rendered, while the LLM retained both the solving and the writing of the solution. This view resonates with recent proposals to distinguish "dependent" from "autonomous" offloading, that is, delegating core thinking from using the tool as a scaffold while retaining agency (Zhu et al., 2026), and to place the role of genAI on a continuum between scaffold and substitute (Cho, 2026). It also resonates with the distinction between executive and instrumental help seeking introduced in Section 2.3 (Nelson-Le Gall, 1981; see also Karabenick & Berger, 2013). "Solve this problem" is executive in this sense, whereas "Explain this term" and "Is my idea correct?" are instrumental, and the two coexisted within single chatlogs and alternated from one task to the next. Thus, we maintain that asking whether a student offloaded their mathematical thinking to an LLM is too coarse to answer with a yes or no. It needs to be asked of each space, and within each space, of each part of the work.

The regulatory lens adds a further qualification. Within our perspective, offloading is not necessarily the student's aim but a means of achieving a self-imposed goal. In Sections 3.1 and 3.3, we argued that the goal of homework should not be presumed (Boekaerts, 1999). A student who set the goal of understanding why a matrix product can vanish, probed the question with the LLM, and then let it write the solution, may have achieved their own goal while offloading the work that the task was designed to elicit. Such use need not be a detour from learning. Studying worked-out solutions is an established route into mathematics, particularly for novices (Atkinson et al., 2000; Kalyuga et al., 2003; Renkl, 2014), and Koichu's (2018) first choice explicitly includes comprehending a worked-out example (Section 3.2). Yet the benefit of an example depends on what learners do with it, most notably on whether they explain it to themselves (Chi et al., 1989).

Thus, an LLM-generated solution may serve a student as a worked example or as a ready-made submission. The chatlogs show which of the two they made of it more clearly than the submission does, although even they do not record what a student did with a solution. Accordingly, a given instance of offloading may be problematic against the goal of the task designer and unproblematic against the goal of the student.

Finally, writing the solution was the part that the students handed over most consistently, keeping, at most, its editing (Section 5.3). As newcomers to linear algebra, they may have regarded an LLM-generated product as meeting the standards better than their own writing, so that handing it over served the mark rather than learning, a proposal that our data do not allow us to test (see Section 6.4).

### 6.2 Task design and the affluence of choice

The second contribution of this study concerns task design, and it rests on a finding that we hold more tentatively than the first. The variation of students' prompting across the three tasks (Section 5.2) is consistent with the event view of regulation adopted in Section 3.1, and it suggests that task features may shape how students divide the work with LLMs. However, the three tasks were given in a fixed order and differed in their mathematical content and in the product they required, so our design cannot separate task design from these other differences (see Section 6.4). Accordingly, we offer what follows as hypotheses for design research rather than as conclusions. In Sections 3.2 and 3.3, we described choice affluence as a matter of what an environment extends to students, and we argued that the instructions accompanying an assignment bound what a student may reasonably put to the LLM. The three tasks illustrate this argument. In Task 1, where the students led the work with minimal guidance, "Solve this problem" appeared in 84% of the chatlogs and checking in 68%, and three of the four abandonments occurred. In Task 2, where the students received a pool of suggested prompts, the exchanges were the longest and most varied, and both major revisions occurred. In Task 3, where students made sense of a given proof, exchanges were short, checking all but disappearing, and no submission departed from the LLM's product beyond minor revisions.

Two hypotheses for task design follow from this comparison. First, supplying prompts that model divisions of labor other than handing over the solving, as in Task 2, may widen the range of choices that students take up. Second, a task that casts the LLM as an explainer does not, by itself, keep the explaining with the student. In Task 3, the product that the students were asked to submit was itself an explanation of a proof, so the LLM's explanations could be carried into the submission almost as they were. Copy-paste was most frequent in this task, although the differences across tasks were small. Accordingly, the relationship between the role a task assigns to the LLM and the product the task requires of the student seems to us worth examining. Where the two coincide, the LLM's product is ready for submission; where they differ, for example, when students are asked to consult an LLM about a conjecture of their own and then to write a proof of it, the student is required to transform what the LLM returns. The inquiries in Section 5.1 point in the same direction, suggesting that tasks inviting students to explore the mathematics around a problem before solving it may be worth designing for.

### 6.3 What chatlogs record of regulation

The third contribution of this study pertains to theory and methodology. Research on self-regulated learning has long struggled with the fact that regulation is ordinarily covert, so researchers have had to reconstruct it from self-reports, think-alouds, or traces left in digital environments (Cleary et al., 2012; Rovers et al., 2019; Winne, 2010). Chatlogs offer such a trace for work with LLMs (Section 3.3). Each prompt inscribes a decision about the task, the roles, and, occasionally, the toolkit and the standards. Our analysis suggests that the constructs of Section 3 are discernible in this record, although unevenly. Every prompt showed task definition and the division of roles, whereas the toolkit and the standards were seldom set in advance and were usually made explicit during monitoring, when an LLM product was compared against a standard that had previously been tacit (Panadero, 2017; Winne & Hadwin, 1998). Indeed, the same holds for judgments of the LLM's mathematics. Only four prompts rejected an LLM reply outright, which echoes reports that students rarely validate what an LLM returns (Bulut & Borromeo Ferri, 2026; Noster et al., 2024). Yet our data suggest that judgments of this kind were voiced in other forms, most often as requests for clarification, and that a study that counts only explicit rejections may underestimate how often students judge what the LLM produced.

The chatlogs are equally telling in what they do not record. Although students could direct the LLM toward other agents (Section 3.3), no prompt did so, even when the LLM's replies were beyond the students' reach (Section 5.2). The students may well have turned to the textbook, the lecture notes, or a classmate outside the exchange, but the chatlogs leave no trace of such a choice.

### 6.4 Limitations and directions for further research

Several limitations should be acknowledged. The study took place in a single course at a teacher education college, with two dozen students. The students were prospective secondary mathematics teachers, so their prompting may reflect a professional interest in explanation that other undergraduates do not necessarily share. The chatlogs were a required part of the submission, so they were themselves products that the students knew their teacher would read, which may have shaped students' prompting. That "Solve this problem" appeared in most chatlogs suggests that the chatlogs were not merely curated to present a learning-oriented use of the LLM, but we cannot rule out such curation. Moreover, our inferences about regulation are bounded by what the chatlogs record (see Section 4.3). The three tasks were also given in a fixed order and differed in their mathematical content, so the variation we attributed to task design may partly reflect when in the semester each task was given and what it was about. The right-to-left direction of Hebrew also scrambled the mathematics in some of the LLM's replies, so that some minor revisions were corrections of the rendering rather than of the mathematics. Finally, the relationship between prompting and submission is correlational. Students who understood the mathematics better may have been both more inclined to share their ideas with the LLM and more motivated to revise what it returned. In other words, students' prompting may have reflected their mathematical understanding rather than shaped it.

These limitations point to directions for further research. First, studies that combine chatlogs with think-alouds or stimulated-recall interviews could reach the regulation that takes place between prompts, including the goals that students set for their work with LLMs. Second, the finding that students most consistently handed over the writing of the solution, submitting LLM products with at most minor revisions, calls for research on how students regard the authorship of a solution when an LLM is available, and on whether writing the solution oneself makes a difference to what is learned. Third, the variation of prompting structures across tasks invites design studies that vary the role a task assigns to the LLM and the product it requires of the student, while counterbalancing task order and holding the mathematical content constant, across courses and over longer periods.

We close by returning to the concern that motivated this study. Horrocks and Shearman (2025) describe the arrival of genAI as an "existential crisis" for mathematics and statistics education. Our study does not dispel the concern behind this description, since most of the solutions in our data stayed close to what an LLM had produced. However, it suggests that the crisis may not be located where it is usually sought. A submitted solution records the product of a student's work with an LLM, not the regulation that produced it. The chatlogs behind the submissions in our data show students defining their own mathematics, asking for explanations, putting conjectures to the LLM, and, in fewer chatlogs than the submissions alone would suggest, relying mainly on the LLM to solve the problem. Accordingly, we propose that research on LLMs in mathematics education shift its question from whether students offload their thinking to which parts of the work they hand over, in which regulatory space, and toward which goal. The answer to these questions, our data indicate, is worked out prompt by prompt, and it cannot be read from the submission alone.

## Statements and Declarations

**Funding.** This study was fully supported by a grant from the Research Grants Council of the Hong Kong Special Administrative Region, China (Project No. 17602126).

**Competing interests.** The authors have no relevant financial or non-financial interests to disclose.

**Ethics approval.** The study was conducted in accordance with the ethical standards of Orot Israel College and with the 1964 Declaration of Helsinki and its later amendments. The research procedures were approved by the college's research ethics committee.

**Consent to participate.** Informed consent was obtained from all individual participants included in the study. Participation was voluntary and separate from the course requirements, and the decision to participate or decline had no bearing on students' grades.

**Consent for publication.** Participants consented to the publication of anonymized excerpts of their chatlogs and submitted solutions. All names reported in the article are pseudonyms, and no identifying information is included.

**Data availability.** The data are not publicly available due to ethical restrictions.

**Author contributions.** Igor' Kontorovich: Conceptualization, Methodology, Formal analysis, Writing – original draft, Writing – review and editing. Miriam N. Wallach: Conceptualization, Investigation, Data curation, Formal analysis, Writing – review and editing. Both authors read and approved the final manuscript.

**Use of generative AI.** During the preparation of this manuscript, the authors used Claude (Anthropic) to edit the English language, proofread the text, and format references. The authors reviewed and revised all output and take full responsibility for the content of the publication.